\documentclass[12pt]{article}
\usepackage{amsthm,amsfonts,amsmath,tikz}
\newtheorem{thm}{Theorem}
\newtheorem{lem}[thm]{Lemma}

\newtheorem*{conj}{Budzik's Conjecture}
\newtheorem*{defn}{Definition}
\begin{document}
\newcommand{\wwd}{\frac{dx_1}{x_1}\wedge\cdots\wedge\frac{dy_\ell}{y_\ell}}
\title{Invariants of Superalgebras as Complex
Integrals}
\author{Allan Berele}
\abstract{In \cite{B} we defined integrals that
approximated the Poincar\'e series of the invariants
and concomiants of the general linear Lie supergroup
or superalgebra.  Budzik suggested in~\cite{Bu} a
way to adapt this method to get the exact Poincar\'e
series.  The purpose of this paper is to prove that
Budzik's ideas are correct.  As a consequence we prove
that the Poincar\'e series are rational functions.}
\maketitle
\section{Budzik's Conjecture}
The basic objects that form the subject of this
paper can be described in several different ways.
The proof of our main theorem will mosty focus on
their definition based on 
the character theory of the symmetic group, but
in order to provide motivation and context, we
also mention some of the connections to trace identities
and to invariant theory.

We will use greek letters $\lambda,$ $\mu$ and $\nu$ to
denote partitions of some unspecified, common integer $n$.  When we
specify the parts of a partition $\lambda=(\lambda_1,\lambda_2,\ldots)$, we
will always intend the parts to be in non-increasing order $\lambda_1\ge
\lambda_2\ge\cdots$.  For each partition $\lambda$ of $n$,
 we let $\chi^\lambda$ be the corresponding character of the symmetric
group $S_n$.  Given $\chi^\mu$ and $\chi^\nu$ the inner tensor product
$\chi^\mu\otimes\chi^\nu$ can be decomposed into a sum of irreducible
parts, whose multiplicities we denoted by $\gamma_{\mu,\nu}^\lambda$,
so that $\chi^\mu\otimes\chi^\nu=\sum \gamma_{\mu,\nu}^\lambda \chi^\lambda$.

Let $k$ and $\ell$ be two fixed non-negative integers,
at least one of which is positive.  Then $H(k,\ell),$ the
$k\times \ell$ hook, is defined to be the set of partitions
$\lambda$ with at most $k$ parts greater than or
equal to $\lambda$, i.e., such that $\lambda_{k+1}
\le\ell$.  The basic objects we will be studying are the
multiplicities occuring are the tensor products
$\sum_{\mu\in H(k,\ell)}\chi^\mu\otimes\chi^\mu$.
The multiplicity $m_\lambda$ of $\chi^\lambda$ in
such a sum is given by
\begin{equation}\label{eq:0}m_\lambda=\sum_{\mu\in
H(k,\ell)}\gamma_{\mu,\mu}^\lambda.\end{equation}
Sometimes we will write $m_\lambda(k,\ell)$ if we need
to specify the $k$ and $\ell$.
We associate to $m_\lambda$ the Poincar\'e series
$$P(k,\ell;n)=\sum m_\lambda s_\lambda(t_1,\ldots,
t_n),$$ where $s_\lambda$ is the Schur function, and
more generally, 
$$P(k,\ell;n,m)=\sum m_\lambda HS_\lambda(t_1,
\ldots,t_n;u_1,\ldots,u_m),$$
where $HS_\lambda$ is the hook Schur function, defined below
for the reader not familiar with them.

The connection to trace identities involves the Kemer
algebras.  Let $E$ be an infinite dimensional Grassmann
algebra in characteristic zero. This algebra has
a natural $\mathbb{Z}_2$-grading, $E=E_0\oplus
E_1$ with respect to which the degree~0 elements
are central and the degree~1 elements anticommute
with each other. The Kemer algebra, $M_{k,\ell}$,
is defined to be the set of $(k+\ell)\times(k+\ell)$
matrices with entries from $E$ of the form $\left(
\begin{smallmatrix}A&B\\C&D\end{smallmatrix}
\right)$ such that $A$ is a $k\times k$ matrix with
entries from~$E_0$, $D$ is an $\ell\times\ell$ matrix
with entries from $E_0$, and $B$ and $C$ have 
entries from $E_1$.  This algebra has a natural
trace function $tr:M_{k,\ell}\rightarrow E_0$ which
takes a matrix as above to $tr(A)-tr(D)$.

To construct generic element for $M_{k,\ell}$ we use the
free supercommutative algebra on the degree~0
variables $x_{ij}^{(\alpha)}$ and the degree~1
variables $e_{ij}^{(\alpha)}$.  Then $X_ \alpha$ is defined
to be the matrix  $\left(
\begin{smallmatrix}A&B\\C&D\end{smallmatrix}
\right)$ in which the dimensions of $A,B,C,D$ are
as above, and the $(i,j)$ entries of $A$ and $D$ are 
$x_{ij}^{(\alpha)},$ and the $(i,j)$ entries of 
$B$ and $D$ are $e_{ij}^{(\alpha)}$.  One next
constructs $F[X_1,\ldots,X_n]$, the algebra generated
by $n$ such matrices and $TR(k,\ell;n)$, the algebra
generated by the traces of the elements of $F[X_1,
\ldots,X_n]$.  The algebra $TR(k,\ell;n)$
has a natural $n$-fold grading and so a Poincar\'e
series in $n$ variables.  It turns out that this 
Poincar\'e series is none other than $P(k,\ell;n)$.

This Poincar\'e series is also the Poincar\'e series for
the invariants of the general linear Lie superalgebra or
supergroup.  In deference to~\cite{Bu} we focus on
 the supergroup, $U(k|\ell)$ (also denoted $PL(k,\ell)$), the group of invertible
 matrices in $M_{k,\ell}$.  This group acts on $M_{k,\ell}$
 by conjugation and, by extension, on the ring of $E_0$-linear
 maps $M_{k,\ell}^n\rightarrow E$.  Such a function $F$ is said
 to be invariant if 
 $$F(ga_1g^{-1},\ldots,ga_ng^{-1})=F(a_1,\ldots,a_n),$$
 for all $a_1,\ldots,a_n\in M_{k,\ell}$ and all $g\in U(k|\ell)$.
 The algebra of invariant functions on $M_{k,\ell}^n,$ polynomial
 in the entries, also has an $n$-fold grading and also has Poincar\'e
 series $P(k,\ell,n)$, the same as above.  In fact it is generated by the trace polynomials.
 
 The series $P(k,\ell;n,m)$ also has a natural interpretation, see~\cite{B05}.  The algebra of $(k+\ell)\times(k+\ell)$
 matrices with entries in~$E$ has a $\mathbb{Z}_2$-grading in which
 $M_{k,\ell}$ is the degree~0 part.  The degree~1
 part would consist of matrices $\left(\begin{smallmatrix}A&B\\C&D\end{smallmatrix}\right)$,with dimensions as
 above, but in which $A$ and $D$ have degree~1
 entries and $B$ and $C$ have degree~0 ones.  This
 algebra has a supertrace function, defined on a
 homogeneous matrix as above as $tr(A)\mp tr(D)$,
 depending on whether the matrix is degree~0 or~1.
 Letting $F[X_1,\ldots,X_n,Y_1,\ldots,Y_m]$ be the
 algebra generated by $n$ generic degree~0 matrices
 and $m$ generic degree~1 matrices, we can then let
 $STR(k,\ell;n,m)$ be the algebra generated by the
 supertraces.  This algebra has an $(n+m)$-fold grading
 and its Poincar\'e series is $P(k,\ell;n,m)$.
 
   We also mention
 that we learned from~\cite{Bu} that these series have signifigance
 in physics: ``These arise in physics e.g. as free partition functions of $U(N|M)$
 gauge theories or supersymmetric indices.''
 
 In \cite{B} we studied these Poincar\'e series using the
 complex integrals.  Let $\Delta$ (not the same as the $\Delta$
 from~\cite{B}) be defined as
$$\Delta=\frac{\prod_{1\le i\ne j\le k}(1-x_ix_j^{-1})
\prod_{1\le i\ne j\le\ell}(1-y_iy_j^{-1})}{\prod_{i=1}^k
\prod_{j=1}^\ell (1+x_iy_j^{-1})(1+x_i^{-1}y_j)}.$$
Let $X=\{x_1,\ldots,x_k\},\ Y=\{y_1,\ldots,y_k\}$ and let $Z_0$ be the set of variables $XX^{-1}\cup YY^{-1}$
and let $Z_1$ be the set of variables $XY^{-1}\cup X^{-1}Y$.
Then for each $\lambda$ we define $m_\lambda'$ or
$m_\lambda'(k,\ell)$ via
\begin{equation}\label{eq:1a}m_\lambda'=\frac1{k!\ell!(2\pi i)^{k+\ell}}
\oint_{|y_i|=1}\oint_{|x_i|=1+\epsilon}
HS_\lambda(Z_0;Z_1)\Delta \frac{dx_1}{x_1}\wedge
\cdots \wedge\frac{dy_\ell}{y_\ell},\end{equation}
where $HS_\lambda$ is the hook Schur function, and the
$\epsilon>0$ is to make sure that all of the poles at $x_i=y_j$
lie inside of the inner integrals.  In~\cite{B} we proved that
$m_\lambda=m_\lambda'$ for all $\lambda$ occuring
with non-zero multiplicity in $\sum_{\mu\in H(k,\ell)}
\chi^\mu\otimes\chi^\mu$,  but not in $\sum_{\mu\in H(k',\ell,)}
\chi^\mu\otimes\chi^\mu,$ for any smaller hook
$H(k',\ell')$.  Based on this, Budzik conjectured in~\cite{Bu}
(slightly reformulated here) that
\begin{conj}$m_\lambda'(k,\ell)=m_\lambda(k,\ell)-m_\lambda(k-1,\ell-1).$\end{conj}
By induction, this would imply that
\begin{equation}\label{eq:A} m_\lambda
(k,\ell)=\sum_{i=0}^{\min(k,\ell)} m_\lambda'(k-i,\ell-i).
\end{equation}
Likewise, since $P'(k,\ell;n)$ is equal to $\sum m_\lambda' S_\lambda
(t_1,\ldots,t_n)$, Budzik's conjecture would imply that $P'(k,\ell;n)=
P(k,\ell;n)-P(k-1,\ell-1;n)$ and so
\begin{equation}\label{eq:B}P(k,\ell;n)=
\sum_{i=0}^{\min(k,\ell)} P'(k-i,\ell-i;n)\end{equation}
More generally, since $P'(k,\ell;n,m)$ is equal to
$\sum m_\lambda' HS_\lambda(t_1,\ldots,t_n;u_1,\ldots,u_m)$ we would
have $P'(k,\ell;n,m)=P(k,\ell;n,m)-P(k-1,\ell-1;n,m)$ and
\begin{equation}\label{eq:C}P(k,\ell;n,m)=
\sum_{i=0}^{\min(k,\ell)} P'(k-i,\ell-i;n,m).\end{equation}
 Our main purpose in this paper is to prove Budzik's conjecture
 and, hence, all of the above equations.
 
 Note that (\ref{eq:1a}) together with Budzik's conjection
 in the form of (\ref{eq:A})
 allows the computation of any $m_\lambda$.  Likewise,
 $P(k,\ell;n)$ and $P(k,\ell;n,m)$ can be computed using
 equations~(\ref{eq:B}) and~(\ref{eq:C}) using Theorem~3.5 from \cite{B}:
 \begin{thm}\label{thm:1} $P'(k,\ell;n,m)$ equals $\frac1{k!\ell!(2\pi
 i)^{k+\ell}}$ times the integral
 $$\oint_{|y_i|=1}\oint_{|x_i|=1+\epsilon}
 \frac{\prod_{z\in Z_0}\prod_i(1+zu_i)\prod_i\prod_{z\in Z_1}
 (1+zt_i)}{\prod_i\prod_{z\in Z_0}(1-zt_i)\prod_i\prod_{z\in Z_1,
}(1-zu_i)}\Delta\wwd$$
and, taking $Y=\emptyset$, 
$P'(k,\ell;n)$ equals $\frac1{k!\ell!(2\pi
 i)^{k+\ell}}$ times the integral
 $$\oint_{|y_i|=1}\oint_{|x_i|=1+\epsilon}
 \frac{\prod_i\prod_{z\in Z_1}
 (1+zt_i)}{\prod_i\prod_{z\in Z_0}(1-zt_i)}\Delta\wwd$$
 \label{thm:1}\end{thm}.

\section{Properties of Hook Schur Functions}

We fix two sets of variables $X=\{x_1,\ldots,x_k\}$ and $Y=\{y_1,\ldots,y_k\}$.
For each partition $\lambda$ we let $HS_\lambda(X;Y)$ be the
hook Schur function on $\lambda$.  This is commonly defined as
$\sum_{\mu\subseteq\lambda}s_\mu(X)s_{\lambda/\mu}(Y)$, where
$s$ is the ordinary Schur function.  Some authors study the super
Schur funcitons
$s_\lambda(X/Y)$ instead of the hook Schur functions $HS_\lambda(X;Y)$.  Each has its
advantages, but they are essentially equivalent, since $HS_\lambda
(X;-Y)=s_\lambda(X/Y)$.

We will need these  properties
of hook Schur funtions:
\begin{thm}[The Hook Theorem] $HS_\lambda(X;Y)=0$ unless $\lambda\in H(k,\ell).$
\end{thm}
Theorem 3.26 of \cite{BR} implies that if $\mu,\nu\in H(k,\ell)$  then $\gamma_{\mu,\nu}^\lambda=0$,
and so $m_\lambda=0$ unless $\lambda\in H(k^2+\ell^2,2k\ell)$.
This means that $P(k,\ell;k^2+\ell^2;2k\ell)$ encodes complete
information about the $m_\lambda$, and no  $P(k,\ell;n,m)$
does for smaller $n,m$.

 A partition in $H(k,\ell)$ with $\lambda_k\ge\ell$ will not 
belong to $H(k-1,\ell-1)$.  Such a partition is called \emph{typical\/} and the
set of such is denoted $H'(k,\ell)$.

The second theorem is due to J\'{o}zefiak and Pragacz, \cite{J}. It
can also be found in \cite{M}, 1.3 Ex.~24.  To 
each partition $\lambda$ we say $(i,j)\in \lambda$ if $\lambda_i\le j$,
or, more geometrically,  if the Young tableau associated to $\lambda$ has a box in row
$i$, column~$j$.  Then $f_\lambda$ is defined to be the product
$\prod (x_i+y_j)$ over all $(i,j)\in \lambda$, with the convention
that $x_i=0$ and $y_j=0$ for $i>k$ or $j>\ell$.  For example, if 
$(k,\ell)=(2,1)$ and $\lambda=(2,1)$, then $f_\lambda=(x_1+y_1)x_1
(x_2+y)$.  J\'{o}zefiak and Pragacz's theorem gives an
alternate description of hook Schur functions.
\begin{thm}[J\'{o}zefiak and Pragacz] $HS_\lambda(X;Y)=$
$$\sum_{\sigma\in S_k\times S_\ell}\sigma\left[
\frac{f_\lambda(X;Y)}{\prod_{i<j}(1-x_j
x_{i}^{-1} )\prod_{i<j}(1-y_{j}
y_{i}^{-1})}\right]$$
\end{thm}
So, for the example given above $HS_{(2,1)}(x_1,x_2;y_1)=$
$$\frac{(x_1+y_1)x_1(x_2+y_1)}{1-x_2x_1^{-1}}+
\frac{(x_2+y_1)x_2(x_1+y_1)}{1-x_1x_2^{-1}}=(x_1+y_1)(x_2+y_1)(x_1+x_2).$$

A weaker version of J\'{o}zefiak and Pragacz's theorem is
called the factorization theorem, and it is due to Berele
and Regev from~\cite{BR}.
\begin{thm}[Berele and Regev] Let $\lambda$ be a typical partition, and define the
partitions $\mu$ and $\nu$, via $\mu=(\lambda_1-\ell,
\ldots,\lambda_k-\ell)$ and $\nu=(\lambda_1'-k,\ldots,
\lambda_\ell'-k)$, see Figure~\ref{fig:1}.  Then $HS_\lambda(X;Y)=\prod_{i,j}(x_i+
y_j)s_\mu(X)s_\nu(Y)$.
\end{thm}
The next theorem is due to Rosas in~\cite{R}.  It gives a
hook Schur function interpretation of 
the coefficients $\gamma_{\mu\nu}^\lambda$ occuring in the
tensor product of the symmetric group characters.
\begin{thm}[Rosas] Given 4 sets of variables $X$, $Y$, $T$ and
$U$,
$$HS_\lambda(XT,YU;XU,YT)=\sum_{\mu,\nu}
\gamma_{\mu,\nu}^\lambda HS_\mu(X;Y)HS_\nu(T;U).$$
In the case of $U=\emptyset$ this says that
$$HS_\lambda(XT;YT)=\sum_{\mu,\nu}
\gamma_{\mu,\nu}^\lambda HS_\mu(X;Y)s_\nu(T).$$
\end{thm}
Finally, we will need this hook generalization of Cauchy's identity
found in~\cite{Rem}
\begin{thm}[Berele and Remmel] Given 4 sets of variables $X$, $Y$, $T$ and $U$,
$$\sum_\lambda HS_\lambda(X;Y)HS_\lambda(T;U)=
\frac{\prod(1+x_iu_j)\prod(1+y_it_j)}{\prod(1-x_it_j)\prod(1-y_i
u_j)},$$
where the $x_i$ run over the elements of $X$, etc.\label{rem}
\end{thm}
\section{The Inner Product}
There is an inner product on symmetric functions in $n$ variables with respect
to which the Schur functions form an orthonormal basis.  It can
be defined using complex integrals via $\langle f,g\rangle=(2\pi i)^{-n}(n!)^{-1}$ times
$$\oint_{|z_1|=1}\cdots\oint_{|z_n|=1} f(z_1,\ldots,z_n)g(z_1^{-1},
\ldots,z_n^{-1})\prod_{i\ne j}(1-z_iz_j^{-1})\frac{dz_1}{z_1}
\wedge\cdots\wedge \frac{dz_n}{z_n}.$$

We now define a new inner product, related to but not the same
as the inner product from~\cite{B}.    In this setting, $f$ and $g$ will be 
functions of two sets of variables $X$ and $Y$, with $|X|=k$ and
$|Y|=\ell$.  This new inner product will be defined via
$\langle f,g\rangle=(2\pi i)^{-(k+\ell)}(k!)^{-1}(\ell!)^{-1}$ times the
integral
$$\oint_{|y|=1}\oint_{|x|=1+\epsilon}f(X;Y)g(X^{-1};Y^{-1})\Delta\wwd.$$
Under this product, the hook Schur functions are not orthonormal.
We proved in~\cite{B} as an easy consequence of the factorization
theorem, that if $\lambda$ and $\mu$ are typical
then $\langle HS_\mu,HS_\nu\rangle =\delta_{\mu,\nu}$.
The main result of this section is the following.
\begin{lem} If $\mu$ and $\nu$ are not both typical, then
$\langle HS_\mu,HS_\nu\rangle=0$.\label{lem:6}
\end{lem}
If $\mu$ or $\nu$ are not in $H(k,\ell)$, then $H_\mu(X;Y)$
or $H_\nu(X;Y)$ will be zero, so we will focus on $\mu,\nu\in 
H(k,\ell)$.
In order to prove this lemma, we need the folowing general fact about
integrals of rational functions.  Just as for polynomials and rational functions, the degree of
a Laurent polynomial is defined to be the degree 
of the highest degree term, and the degree of
a quotient of Laurent polynomials is the difference
between the degrees of the numerator and 
denominator.
\begin{lem} Let $f=f(z_1,\ldots,z_n)$ and $g=g(z_1,\ldots,z_n)$
be Laurant polynomials with leading terms of degrees $a$ and $b$,
respectively, with $a<b-1$.  Let $|z|=r$ be a circle containing all of the roots
of $g(z)$ in its interior.  Then $\oint_{|z|=r} f(z)/g(z)dz=0$.
\label{lem:7}\end{lem}
\begin{proof} By Cauchy's residue theorem $\oint_{|z|=r} f(z)/g(z)dz=0$ is not changed if we replace the circle $|z|=r$ with
a larger circle $|z|=R$.  For large $R$ when $|z|=R$, $|f(z)/g(z)|$ is 
$O(R^{a-b})$.  Since the circumfrence of this circle is $2\pi R$,
the absolute value of the integral is $O(R^{a-b+1})$, which tends to
zero.
\end{proof}
We now analyze the fraction
\begin{equation}\label{eq:1}
\frac{HS_\mu(X;Y)HS_\nu(X^{-1};Y^{-1})\prod(1-x_ix_j^{-1})
\prod(1-y_iy_j^{-1})}{\prod(1-x_iy_j^{-1})(1-x_i^{-1}y_j)\prod x_i
\prod y_j}\end{equation}
\begin{lem} If $\mu$ is not typical, then the integrand (\ref{eq:1}) can be written as a sum of terms of the form $x^\alpha y^\beta\prod(x_i+y_j)^{-2}$,
each of degree at most $-2$ in one of the $x_i$.
\end{lem}
\begin{proof}
Isolating the term  
\begin{equation}\label{eq:1b}\frac{\prod(1-x_ix_j^{-1})
\prod(1-y_iy_j^{-1})}{\prod(1-x_iy_j^{-1})(1-x_i^{-1}y_j)\prod x_i
\prod y_j},\end{equation}
the numerator  is of degree $k-1$ in each $x_i$, and the denominator is of degree $\ell+1$, so the degree of the fraction is
$k-\ell-2$ in each $x_i$.

By J\'ozefiak and Pragacz's Theorem, $HS_\mu(X;Y)$ is a sum of the terms
\begin{equation}\label{eq:1c}\sigma\left[f_\mu(X,Y)\prod_{i<j}(1-x_j
x_{i}^{-1} )^{-1}\prod_{i<j}(1-y_j
y_{i}^{-1})^{-1}\right].\end{equation}
Since $\mu$ is not typical,  $\mu_k\le \ell-1$, so $f_\mu(X;Y)$ is of
degree at most $\ell-1$ in $x_k$; also $\prod_{i<j}(1-x_{j}
x_{i}^{-1} )^{-1}$ will be of degree $k-1$ in $x_k$, so (\ref{eq:1c}) will
be of degree at most $\ell-k$ in $x_k$.  Under a permutation $\sigma$, $f_\lambda(X;Y)/\prod(1-x_jx_i^{-1})$ will be of degree at most $\ell-k$
  in $x_{\sigma(k)}$, and the product with (\ref{eq:1b}) at
most~$-2$.  
\end{proof}
\begin{proof}[Proof of Lemma~\ref{lem:6}]
By the previous lemma, if $\mu$ is not typical, then
$\langle HS_\mu(X;Y)HS_\nu(X;Y)\rangle$ is a linear 
combination of integrals of the form
$$\oint_{|y|=1}\oint_{|x|=1+\epsilon}\frac{x^\alpha y^\beta}{\prod(x_i+y_j)^2}\wwd,$$
such that for at least one $i$, $x_i^{\alpha_i}\prod_j(x_i+y_j)^{-2} $the degree
in $x_i$ is less than or equal to $-2$.  The integral with respect to the $x_j$ can be broken up into the product
of integrals
$$\oint\frac{x_j^{\alpha_j}}{\prod_m (x_j+y_m)^2}\frac
{dx_j}{x_j}.$$
For $j=i$, the fraction will be of degree at most
$-2$ and so must be zero by Lemma~\ref{lem:7}.
\end{proof}

\section{Proof of Budzik's Conjecture}
The proof is now easy:
\begin{align*}
m_\lambda(k,\ell)&-m_\lambda(k-1,\ell-1)\\&=\sum_\mu\{\gamma_{\mu,\mu}|\mu\in H(k,\ell)/H(k-1,\ell-1)\}\\
&=\sum_{\mu,\nu}\gamma_{\mu,\nu}^\lambda
\langle HS_\mu(X;Y),HS_\nu(X;Y)\rangle\mbox{ (by Lemma \ref{lem:6})}\\
&=\langle \sum_{\mu,\nu}\gamma_{\mu,\nu}^\lambda HS_\mu(X;Y)HS_\nu(X^{-1};Y^{-1}),1\rangle 
\\
&=\langle HS_\lambda(XX^{-1}\cup YY^{-1};XY^{-1}\cup X^{-1}Y),
1\rangle\mbox{ (by Rosas's Theorem)}\\
&=\langle HS_\lambda(Z_0;Z_1),1\rangle\\
&=m_\lambda'
\end{align*}
Note that using Theorem~\ref{rem} we can directly
compute the formula for $P'(k,\ell;n,m)$:
\begin{align*}
P'(k,\ell;n,m)&=\sum m_\lambda HS_\lambda(T,U)\\
&=\sum \langle HS_\lambda (Z_0;Z_1),1\rangle
HS_\lambda(T;U)\\
&=\langle \sum HS_\lambda(Z_0;Z_1)HS_\lambda(T;U)
,1\rangle
\end{align*}
The sum is the integrand of Theorem~\label{thm:1}
by the Berele-Remmel Theorem, Theorem~\ref{rem}.
\section{Concomitants}
In \cite{B} we studied not only the invariants of $U(k|\ell)$, but
also the concomitants.  Each of the three definitions of $m_\lambda$ have analogues to multiplicities we call 
$\bar{m}_\lambda$.  In terms of tensor products $\sum 
\bar{m}_\lambda$ is the multiplicity of $\chi^\lambda$ in
$\sum_{\mu\in H(k,\ell)}(\chi^\mu\otimes\chi^\mu)\downarrow,$
where the arrow denotes inducing down from each $S_n$ to
$S_{n-1}$.  For the second one, we use $\bar{F}[X_1,\ldots,X_n]$,
defined to be the tensor product $F[X_1,\ldots,X_n]\otimes
TR(k,\ell;n)$.  This algebra is $n$-graded and has Poincar\'e
series $\sum \bar{m}_\lambda s_\lambda(t_1,\ldots,t_n)$,
denoted $\bar{P}(k,\ell;n)$.  
And, for the third one, we use the algebra of $U(k|\ell)$-invariant
functions $M_{k,\ell}^n\rightarrow M_{k,\ell}$.  This turns out
to be isomorphic to the previous, and so also has Poincar\'e
series $\bar{P}(k,\ell;n)$.  Likewise for $\bar{P}(k,\ell;
n,m)$.

There is an analogue of Budzik's conjecture in this case and it,
together with the analogues of equations~(\ref{eq:A}),
(\ref{eq:B}) and~(\ref{eq:C}) turn out to be true.  Let
$Z=Z_0\cup Z_1.$  The 
analogue of $m_\lambda'$ is
$$\bar{m}_\lambda'=\frac1{k!\ell!(2\pi i)^{k+\ell}}
\oint_{|y_i|=1}\oint_{|x_i|=1+\epsilon}
HS_\lambda(Z_0;Z_1)\Delta\sum_{z\in Z} z \frac{dx_1}{x_1}\wedge
\cdots \wedge\frac{dy_\ell}{y_\ell},$$
and the analogue of $P'(k,\ell;n)$ is $\bar{P}'(k,\ell;n)$
which equals $\frac1{k!\ell!(2\pi
 i)^{k+\ell}}$ times the integral
 $$\oint_{|y_i|=1}\oint_{|x_i|=1+\epsilon}
 \frac{\prod_i\prod_{z\in Z_1}
 (1+zt_i)\sum_{z\in Z}z}{\prod_i\prod_{z\in Z_0}(1-zt_i)}\Delta\wwd$$
Note that in both cases, the integrals are gotten from those
of $m_\lambda'$ and $P'(k,\ell;n)$ by adding in a factor
of $\sum_{z\in Z}z$.  $P'(k,\ell;n,m)$ is defined in the same
way, but we don't write it out because of its length.

To start the proof, we relate $\bar{m}_\lambda'$ to $\bar{P}'(k,\ell;n,m)$
and $\bar{P}'(k,\ell;n)$.
\begin{lem}  The functions $\bar{P}'(k,\ell;n,m)$
and $\bar{P}'(k,\ell;n)$ satisfy $\bar{P}'(k,\ell;n)=\sum \bar{m}_\lambda' s_\lambda(T)$ 
and $\bar{P}'(k,\ell;n,m)=\sum m'_\lambda HS_\lambda(T;U)$.
\label{lem:10}\end{lem}
\begin{proof}
For simplicity, we prove only the first statement.  The proof of the second
is essentially the same.  
\begin{align*} 
\sum_\lambda \bar{m}'_\lambda s_\lambda(T)&=
\sum_\lambda \langle HS_\lambda(Z_0;Z_1)\sum z,1\rangle s_\lambda(T)\\
&=\langle HS_\lambda(Z_0;Z_1)s_\lambda(T)\sum z,1\rangle\\
&=\langle \prod_i\prod_{z\in Z_0}(1+zt_i)\prod_i\prod_{z\in Z_1}
(1+zu_i)\times\\ &\quad\prod_i\prod_{z\in Z_1}(1+zt_i)^{-1}\prod_i\prod_{z\in Z_0}
(1+zu_i)^{-1}\sum z,1\rangle\\&=P'(k,\ell;n)
\end{align*}
The third inequality is Theorem~\ref{rem}.
\end{proof}

The next lemma is a standard trick in trace identity theory.
\begin{lem} $\bar{P}(k,\ell;n)$ can be gotten from $P(k,\ell;n+1)$
via 
$$\bar{P}(k,\ell;n)=\left.\frac{\partial P(k,\ell;n+1)}{\partial t_{n+1}}\right|_{t_{n+1}=0}.$$\label{lem:11}
\end{lem}
\begin{proof}
We map $\bar{F}[X_1,\ldots,X_n]\rightarrow TR(k,\ell;n+1)$ via
$f\mapsto tr(X_{k+1}f)$.  Since the trace is non-degenerate, this is
an $n$-grading preserving linear isomorphism onto the part of $TR(k,\ell;n+1)$
of degree 1 in $X_{k+1}$.  Hence, the Poincar\'e series of 
$\bar{F}[X_1,\ldots,X_n]$ is equal to the Poincar\'e series of the subspace
of $TR(k,\ell;n+1)$ of degree~1 in $X_{k+1}$.  The operation of
taking derivative with respect to $t_{k+1}$ and then setting $t_{k+1}$
equal to 0 computes the part of $\bar{P}(k,\ell;n+1)$ of degree~1
in~$t_{k+1}$.
\end{proof}

The next lemma does this computation.
\begin{lem} For each $i$, $\frac{\partial}{\partial t_{i+1}}
P'(k,\ell;i+1)\vert_{t_{i+1}=0}=\bar{P}'(k,\ell;i)$.\label{lem:12}
\end{lem}
\begin{proof} The integrand computing $P'$  from Theorem~\ref{thm:1},
taking $n=i$, is a
product of terms not involving $t_i$ times
$\prod_{z\in Z_1}(1+zt_i)\prod_{z\in Z_0}(1-z_it)^{-1}$.  Taking
the partial derivative with respect to $t_i$ and then setting
$t_i=0$ gives a sum of terms $\sum_{z\in Z_1}z+\sum_{z\in Z_0}z=
\sum_{z\in Z}z.$.
\end{proof}
Here is the main theorem of this section.
\begin{thm} For each $\lambda$, $\bar{m}_\lambda'(k,\ell)=
\bar{m}_\lambda(k,\ell)-\bar{m}_\lambda(k-1,\ell-1)$, and so
$\bar{m}_\lambda(k,\ell)=\sum m_\lambda'(k-i,\ell-i)$.

For each $n$, $\bar{P}'(k,\ell;n)=\bar{P}(k,\ell;n)-\bar{P}(k-1,\ell-1;n)$,
and so $\bar{P}(k,\ell;n)=\sum \bar{P}'(k-i,\ell-i;n)$.  And, more
generally, $\bar{P}'(k,\ell;n,m)=\bar{P}(k,\ell;n,m)-\bar{P}(k-1,\ell-1;n,m)$,
and so $\bar{P}(k,\ell;n,m)=\sum \bar{P}'(k-i,\ell-i;n,m)$.
\end{thm}
\begin{proof} We begin with the case of $\bar{P}(k,\ell;n).$
\begin{align*}
\bar{P}(k,\ell;n)&=\left.\frac{\partial P(k,\ell;n+1)}{\partial t_{n+1}}\right|_{t_{n+1}=0}\mbox{ (by Lemma~\ref{lem:11})}\\
&=\sum \left.\frac{\partial P'(k-i,\ell-i;n+1)}{\partial t_{n+1}}\right|_{t_{n+1}=0}\mbox{ (by (\ref{eq:B}))}\\
&=\sum \bar{P}'(k-i,\ell-i;n) \mbox{ (by Lemma~\ref{lem:12})},
\end{align*}
which implies $\bar{P}'(k,\ell;n)=\bar{P}(k,\ell;n)-\bar{P}(k-1,\ell-1;n)$.

This implies, using Lemma~\ref{lem:10}, that 
$$\sum \bar{m}_\lambda'(k,\ell)
s_\lambda=\sum \bar{m}_\lambda'(k,\ell)s_\lambda
-\sum \bar{m}_\lambda'(k-1,\ell-1)s_\lambda.$$  Since
the Schur functions are linearly independent, this implies that
 $\bar{m}_\lambda'(k,\ell)=\sum \bar{m}_\lambda(k,\ell)-
\sum \bar{m}(k-1,\ell-1)$, which in term implies
$\bar{P}'(k,\ell,n,m)=P(k,\ell;n,m)-P(k-1,\ell-1;n,m)$.

\end{proof}
\section{Algebraic Properties of the Poincar\'e Series}
Our main goal in this section is to prove that each of the
Poincar\'e series studied in this paper, namely,
 $P(k,\ell;n)$, $P(k,\ell;n,m)$, $\bar{P}(k,\ell;
n)$, $\bar{P}(k,\ell;n,m)$, $P'(k,\ell;n)$, $P'(k,\ell;n,m)$
$\bar{P}'(k,\ell;n)$ and $\bar{P}'(k,\ell;n,m)$ is a nice
rational function, by which we mean that each is a rational function
whose denominator can be written as a product of terms
of the form $1-w$, where $w$ is a monic monomial. We will
also prove that the monomials $w$ each have degree at most
$k+\ell$.  Moreover, we prove that if $\ell=0$,
each $w$ has degree at most $\max\{k,
\ell\}$; and in any case $w$ has even degree in the $u_i$.\footnote{The 
computed value of $\bar{P}'(1,1;2,2)$
in \cite{B}
has a factor of $(1-u_1)(1-u_2)$ in the 
denominator.  However, multiplying
numerator and denominator by $(1+u_1)(1+u_2)$ gives a fraction with $(1-u_1^2)
(1-u_2^2)$ in the denominator, both terms of even degree in $u_i$ and both of 
degree $\le 2$.}  Note that the first four in our list (the unprimed ones) are
sums of the last four, so it suffices to prove that these four
are nice rational.  For ease of exposition we will first do the
case of $P'(k,\ell;n)$ and then indicate
the changes needed for the remaining
cases.  

By Theorem~\ref{thm:1}, $P'(k,\ell;n)$ is the integral of
the product of 
\begin{equation}\label{eq:9}
\frac{\prod_{i,j,\alpha}(1+x_i^{-1}y_jt_\alpha)(1+x_iy_j^{-1}
t_\alpha)
\prod_{i\ne j}(1-x_ix_j^{-1})\prod_{i\ne j}(1-y_iy_j^{-1})}
{\prod_{i,j,\alpha}(1-x_ix_j^{-1}t_\alpha)\prod_{i,j,\alpha}(1-y_iy_j^{-1}t_\alpha)\prod_i x_i\prod_j
y_j}\end{equation}
and
\begin{equation}\label{eq:10}\prod_{i,j}(1+x_i^{-1}y_j)^{-1}
(1+x_iy_j^{-1})^{-1}.\end{equation}
We first focus on this latter product. For each $i,j,$  $x_i$ has
norm $1+\epsilon$ and  $y_j$ has smaller norm, 1, so we
can write each $(1+x_i^{-1}y_j)^{-1}(1+x_iy_j^{-1})^{-1}$
as $x_i^{-1}y_j(1+x_i^{-1}y_j)^{-2}=x_i^{-1}y_j\sum_{m=0}^\infty (m+1)(-x_i^{-1}y_j)^m.$  Hence, (\ref{eq:10}) equals
$\prod_i x_i^{-\ell}\prod_j y_j^k \prod_{i,j}\sum_{m=0}^\infty (m+1)
(-x_i^{-1}y_j)^m$.

In (\ref{eq:9}) the degree in each $x_i$ is $n(k+\ell)-1$, so
by Lemma~\ref{lem:7}, the integral of (\ref{eq:9}) times
$(-x^{-1}y_j)^m$ will be zero if 
$m\ge n(k+\ell)+1$.  Letting $N$ equal to $n(k+\ell)+1$,
we can replace each $\sum_{m=0}^\infty (m+1)
(-x_i^{-1}y_j)^m$ by $\sum_{m=0}^N (m+1)
(-x_i^{-1}y_j)^m$ without effecting the integral.  Hence,
$P'(k,\ell;n)$ is a linear combination of (many) integrals with
integrands of the form 
\begin{equation}\label{eq:11}
\frac{x_1^{i_1}\cdots x_k^{i_k}y_1^{j_1}\cdots y_\ell^{j_\ell}
t_1^{\alpha_1}\cdots t_n^{\alpha_n}}{\prod_{i,j,\alpha}(1-x_ix_j^{-1}
t_\alpha)\prod_{i,j,\alpha}
(1-y_i y_j^{-1} t_\alpha)}.\end{equation}
That such an integral is a nice rational function follows directly
from Theorem~3.4 of~\cite{V}, which in turn is based on Stanley's Corollary~3.8 in \cite{S}.  Van Den Bergh's theory associates directed graphs
 to integrals of this form.  In our
case, the graph would have vertices
$v[x_i]$ and $v[y_j]$, corresponding
to the variables $x_1,\ldots,y_\ell$, and
  with an edge $e[x_i,x_j,t_\alpha]$ connecting $v[x_i]$ with $v[x_j]$ for each term $(1-x_i
x_j ^{-1}t_\alpha)$  in
the denominator, and an edge $e[y_i,y_j,
t_\alpha]$ from $v[y_i]$ to $v[y_j]$
 for each $(1-y_iy_j^{-1}t_\alpha)$.  Then his Theorem~3.4 
  not only states that the integral is a nice rational function, but
  identitifies the denominator has a product of terms of the form
  $1-t^C$, where $C$ is an oriented cycle in the graph, and where $t^C$ means the
product of the $t_\alpha$ corresponding
to edges in $C$.  Since a
  cycle cannot have repeated vertices,
and since edges connect only the $v[x_i]$ to each other
 or the $v[y_j]$ to each other,  $|C|\le \max\{k,\ell\}$, and so
  the degree of $t^C$ is at most $\max
\{k,\ell\}$.

The case of $\bar{P}'(k,\ell;n)$ is essentially the same.  As for the cases
of $P'(k,\ell;n,m)$ and $\bar{P}'(k,\ell;n,
m)$, the integrands~(\ref{eq:11}) would
have extra factors of 
$$\prod (1-x_iy_j^{-1}u_\alpha)(1-x_i^{-1}
y_ju_\alpha).$$
In terms of the graph theory, we would
add to the previous graph edges $e[x_i,y_j,u_\alpha]$
 for  each 
$(1-x_iy_j^{-1}u_\alpha)$ and $e[y_i,x_j,u_\alpha]$ for each
$(1-x_i^{-1}y_ju_\alpha)$.  The effect of these additional
edges would be
that cycles could now have length as
great as $k+\ell$, the number of vertices,
and the degree in the edges corresponding
to the $u_\alpha$ must be even.
\section{Computation of $P'(k,\ell;0,1)$}
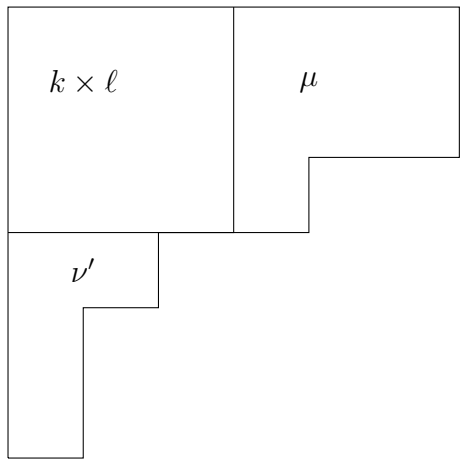
\begin{figure}
\begin{center}
\begin{tikzpicture}
\draw (0,0)--(1,0)--(1,2)--(2,2)--(2,3)--(4,3)--(4,4)--(6,4)--(6,6)--(0,6)--(0,0);
\draw (0,3)--(3,3)--(3,6);
\node at (1,5){$k\times\ell$};
\node at (4,5){$\mu$};
\node at (1,2.5){$\nu'$};
\end{tikzpicture}
\caption{A Typical Partition}\label{fig:1}
\end{center}
\end{figure}

In \cite{Bu} Budik computed $P'(k,\ell,1,0)$, the Poincar\'e series
of the ring generated by the traces of a single generic element
of $M_{k,\ell}$:
\begin{equation}\label{bud}P'(k,\ell;1)=\frac{t^{k\ell}}
{\prod_{i=1}^k (1-t^i)\prod_{j=1}^\ell
 (1-t^j)},\end{equation}
Also computed by Molev and Mukhin in~\cite{Mo}, Orellana and 
Zabrocki in~\cite{O}.  A different formula for $P(k,\ell;n)$ was
given by Sergeev and Veselov in~\cite{Ser}.  We give yet another
derivation based on Theorem~5.2 of~\cite{B}, our proof being completely combinatorial. 
\begin{thm}[=Theorem 5.2 of \cite{B}] $m_{(n)}(k,\ell)$ equals the
number of partitions of~$n$ in the hook $H(k,\ell)$; and $m_{(1^n)}
(k,\ell)$ equals the number of self-conjugate partitions of~$n$ in
$H(k,\ell)$.\label{thm:14}
\end{thm}
It follows that $m'_{(n)}(k,\ell)=m_{(n)}(k,\ell)-m_{(n)}
(k-1,\ell-1)$ equals the number of partitions in $H(k,\ell)$, but
not in $H(k-1,\ell-1)$.  These are the typical partitions, $H'(k,\ell)$.
A Young diagram of a typical partition consists of three parts:  
The $k\times \ell$ rectangle, a partition of height at most $k$ to
the right of it, and the conjugate of a partition of height at most
$\ell$ below it.  See figure~\ref{fig:1}.  Since the generating
function for partitions of height at most $k$ is $\prod_{i=1}^k (1-x^i)^{-1}$, (\ref{bud}) follows.

We now turn our attention to $P'(k,\ell;0,1)$.  We will take $k\ge \ell$.
 This is the
Poincar\'e series of the supertrace ring generated by one
matrix of the form $\left(\begin{smallmatrix}A&B\\C&D
\end{smallmatrix}\right)$, where $A$ is a $k\times k$ matrix,
$D$ is an $\ell\times\ell$ matrix, $A$ and $D$ have generic
odd Grassmann entries, and $B$ and $C$ have entries commuting
indeterminants.  

As above, 
$$m_{(1^n)}'(k,\ell)=m_{(1^n)}(k,\ell)-m_{(1^n)}(k-1,\ell-1),$$
which, by Theorem~\ref{thm:14}, equals the number of self-conjugate partitions of~$n$
in the hook $H(k,\ell)$, but not in the hook $H(k-1,\ell-1)$.  

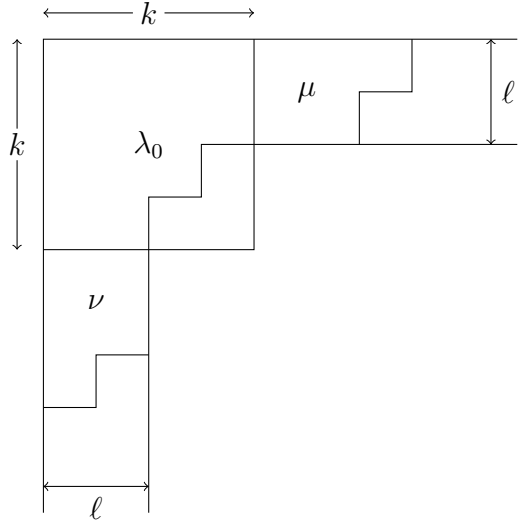
\begin{figure}
\begin{center}
\begin{tikzpicture}[scale=.7]
\draw(0,0)--(0,9)--(9,9);
\draw (0,2)--(1,2)--(1,3)--(2,3);
\draw (2,5)--(2,0);
\draw (0,5)--(4,5)--(4,9);
\draw(4,7)--(9,7);
\draw (6,7)--(6,8)--(7,8)--(7,9);
\node at (2,7){$\lambda_0$};
\node at (1,4){$\nu$};
\node at (5,8){$\mu$};
\draw[<->] (0,.5)--(2,.5);
\node[below] at (1,.5){$\ell$};
\draw[<->]  (8.5,9)--(8.5,7);
\node[right] at (8.5,8){$\ell$};
\draw[<-] (0,9.5)--(1.7,9.5);
\draw[->] (2.3,9.5)--(4,9.5);
\node at (2,9.5){$k$};
\draw[<-](-.5,5)--(-.5,6.7);
\draw[->] (-.5,7.3)--(-.5,9);
\node at (-.5,7){$k$};
\draw (2,5)--(2,6)--(3,6)--(3,7)--(4,7);
\end{tikzpicture}
\caption{Construction of $\lambda_0$, $\mu$ and
$\nu$}\label{fig:2}
\end{center}
\end{figure}

If a self-conjugate partition $\lambda$ lies in the hook $H(k,\ell)$,
it will also lie in $H(\ell,k)$, meaning that $\lambda_{k+1}\le \ell$
and $\lambda_{\ell+1}\le k$.  In terms of Young diagrams, such
a $\lambda$ can be broken into three parts.
\begin{defn}
Given a partition $\lambda$ we define $\lambda_0$ to
be partition $\lambda\cap (k^k)$, namely, the partition 
whose Young diagram is the intersection of the Young diagram
of~$\lambda$ with that of the square $(k^k)$;   $\mu=\mu(\lambda)$
to be the partition to the right of $\lambda_0$; and $\nu=\nu(\lambda)$
to be the conjugate of the partition below it, see Figure~\ref{fig:2}.
\end{defn}

Note that $\lambda$ is self-conjugate if and only if $\lambda_0$
is self-conjugate and $\mu(\lambda)=\nu(\lambda)$.
\begin{lem} The correspondence $\lambda\rightarrow(\lambda_0,
\mu(\lambda))$ is a one-to-one correspondence between
self-conjugate partitions in $H'(k,\ell)$ and pairs of partitions, where
$\mu$ is of height at most $\ell$  and $\lambda_0$ is a self-conjugate partition
 containing $(k^\ell,\ell^{k-\ell})$ and contained in
$(k^k)$.\end{lem}
\begin{proof} If $\lambda\notin H(k-1,\ell-1)$ then $\lambda_k\ge
\ell$ and, since $\lambda$ is self-conjugate, $\lambda_\ell\ge
k$.  The latter implies that the first $\ell$ rows of $\lambda_0$
must equal $k$, and the former that the first $\ell$ columns
also equal $k$.  Hence, $\lambda_0$ contains $(k^\ell,\ell^{k-\ell})$.

The condition $\lambda\in H(k,\ell)\cap H(\ell,k)$ implies that
$\mu$ has height at most $\ell$, and the conditions $\lambda_\ell
\ge k$ and $\lambda_k\ge \ell$ imply that any partition of height at
most $\ell$ is  possible.
\end{proof}

\begin{thm} $P'(k,\ell;0,1)=u^{2k\ell-\ell^2}\prod_{i=1}^{k-\ell} (1+u^{2i-1})\prod_{i=1}^\ell (1-u^{2i})^{-1}$.
\end{thm}
\begin{proof} The generating function for partitions of height at most
$\ell$ is $\prod_{i=1}^\ell (1-u^i)^{-1}$, so the generating function
for equal pairs of them is $\prod_{i=1}^\ell (1-u^{2i})^{-1}$.  As
for $\lambda_0$, $\lambda_0$ is made up of two parts:  
$(k^\ell,\ell^{k-\ell})$ and a self-conjugate partition contained
in the lower right $(k-\ell)\times(k-\ell)$ square.  The former
has degree $k\ell+\ell(k-\ell)=2k\ell-\ell^2$, and the latter has
generating function $\prod_{i=0}^{k-\ell} (1+u^{2i-1}).$
\end{proof}

Budzik related the formula for $P'(k,\ell;n)$ to various important
combinatorial identities.  The same came be done for $P'(k,\ell;0,1)$.
The simplest is based on $k=\ell$, in which case
$$P'(k,k;0,1)=\frac{u^{k^2}}{[u^2]_k},$$
the generating function for self-conjugate partitions in $H'(k,k)$. 
Using $P(k,k;0,1)=\sum_{n=0}^k P'(n,n;0,1)$ and letting
$k$ go to infinity we get the well-known identity
$$\prod_{n=1}^\infty (1+u^{2n-1})=\sum_{k=0}^\infty \frac{u^{k^2}}{[u^2]_k}$$
since $P(\infty,\infty;0,1)$ is the generating function for all self-conjugate partitions.

Similarly, using the identity $P(n+\ell,\ell;0,1)=\sum_{i=0}^\ell
P'(n+i,i;0,1)$ and letting $\ell$ go to infinity, we get
$$\prod_{i=1}^\infty (1+u^{2i-1})=\sum_{i=0}^\infty
\frac{u^{i(2n+i)}\prod_{t=1}^i (1+u^{2t-1})}{[u^2]_i}.$$

\end{document}